\magnification=\magstep1
\input amstex
\documentstyle{amsppt}
\redefine\a{\alpha}
\redefine\b{\beta}
\redefine\g{\gamma} 
\redefine\d{\delta}
\redefine\o{\omega} 
\redefine\s{\sigma}
\redefine\k{\kappa}
\redefine\R{\Bbb R} 
\redefine\C{\text{\rm Club}}
\define\club{\text{\tt pred}}
\define\CL{\Bbb {CL}}

\define\rng{\text{\rm rng}}

\define\cov{\text{\tt cov}}
\define\add{\text{\tt add}}
\redefine\colon{{:}\;}

\pageheight{7.1in}
\baselineskip 15pt

\topmatter
\title
An interpolation theorem
\endtitle
\author
Jind\v rich Zapletal
\endauthor
\affil
University of Florida
\endaffil
\abstract
If $\frak x$ is a tame cardinal invariant and $\frak x=\aleph_1$ implies a certain prediction principle WCG on $\omega_1,$ this happens only because $\frak x=\aleph_1$ implies $\frak b=\aleph_1$ and this in turn implies the principle WCG.
\endabstract
\address
PO Box 118105, Gainesville FL 32611-8105
\endaddress
\thanks
The author is partially supported by grants GA \v CR 201-00-1466, NSF DMS-0071437 and a CLAS UF research award.
\endthanks
\email
zapletal\@math.ufl.edu
\endemail
\endtopmatter


\document

\head{0. Introduction}\endhead

This note is an attempt at quantifying the old belief that the interplay between cardinal invariants of the continuum and the prediction principles at $\aleph_1$ is quite simple. 
I will deal with one of the very few prediction principles I can handle at this point, the weak club guessing principle WCG.

\definition {0.1. Definition}
The {\bf weak club guessing principle} WCG is the statement that there is a collection
of $\aleph_1$ many subsets of $\o_1$ of ordertype $\o$ such that
every club subset of $\o_1$ has an infinite intersection with one of them.
\enddefinition

It is not difficult to observe that $\frak b=\aleph_1$ implies WCG. 
The point of this note
is to assert that this is the only implication of its sort\colon

\proclaim {0.2. Theorem} 
Whenever $\frak x$ is a tame cardinal invariant, if ZFC+LC
proves that $\frak x=\aleph_1$ implies WCG then ZFC+LC proves that $\frak b\leq 
\frak x$. 
In other words, the provable implication $\frak x=\aleph_1\to$ WCG is a
composition of the provable implications $\frak x=\aleph_1\to\frak b=\aleph_1$
and $\frak b=\aleph_1\to$ WCG.
\endproclaim

This is not an unprecedented situation. 
If $\frak x$ is a tame cardinal invariant and ZFC+LC proves that $\frak x=\aleph_1$ implies the failure of the Borel conjecture (BC), then that implication is the composition of two ZFC+LC provable implications $\frak x=\aleph_1\to\frak b=\aleph_1$ and $\frak b=\aleph_1\to\lnot$BC. 
To see this suppose that ZFC+LC does not prove $\frak b\leq \frak x$. 
By \cite {Z} and \cite {L}, within the consistent theory ZFC+LC+$\frak x<\frak b$ it is possible to conclude that in the Laver model $\frak x=\aleph_1$ holds in conjunction with BC. 
Thus ZFC+LC does not prove the implication $\frak x=\aleph_1\to\lnot$BC. 
From this argument it is clear that the way to establish Theorem 0.2 is to prove

\proclaim {0.3. Theorem}
(ZFC+LC) There is a partial order $P$ with the following property.
If $\frak x$ is a tame cardinal invariant and $\frak x<\frak b$ holds in some
forcing extension, then  $\aleph_1=\frak x<\frak b=\aleph_2$ 
holds in the $P$-extension
in conjunction with the failure of the weak club guessing principle.
\endproclaim

It is also clear that the methods of \cite {Z} will have to be adjusted--all the extensions described in that paper use $<\o_1$-proper notions of forcing and therefore preserve even much stronger prediction principles than WCG. 

The main open question in the area is whether similar theorems can be proved for
other prediction principles, most notably the Ostaszewski principle. 
In that
particular case it is natural to conjecture that the interpolating statement should be
$\frak c=\aleph_1$, however the technique used in this note does not immediately
generalize to prove such a thing.

The notation in this paper follows \cite {J}. 
A tame cardinal invariant is one defined as
$\min\{|A|\colon A\subset\R,\phi(A)\land\forall x\in\R\exists y\in A\ \theta(x,y)\}$,
where the quantifiers in the formula $\phi$ range over natural numbers or elements of the set $A$ and $\theta$ is a projective formula making no mention of the set $A$, \cite {Z}. 
An $\epsilon$-number
is a countable limit ordinal closed under ordinal exponentiation. 
Borel sets are constantly confused with their definitions. 
Given a Borel set $B$ and a partial order $P,$ the symbol $\dot B$ denotes the $P$-name for the set of reals with the same Borel definition as $B.$ 
Thus $P\Vdash\check B=\dot B\cap V.$ 
LC denotes a``suitable large cardinal
hypothesis'', in Theorem 0.3 it can be specified to be ``a proper class of measurable
Woodin cardinals'', for Theorem 0.2 ``$\omega_1$ many Woodin cardinals'' is sufficient.

\head {1. The ideals $\club_\b$}\endhead

\definition {1.1. Definition}
Let $\beta$ be a countable limit ordinal. 
By $\C_\beta$ I will denote the collection
of all closed subsets of $\beta$ of ordertype $\b$. 
Identifying sets with their characteristic functions, the collection $\C_\b$ is viewed as a Borel subset of the Polish space $2^\beta$ equipped with the product topology.
The symbol $\club_\b$ denotes the ideal $\s$-generated by the
sets $A_d=\{e\subset\b\colon e\cap d$ is infinite$\}$ as $d$ ranges over all
subsets of $\b$ of ordertype $\o$.
\enddefinition

\proclaim {1.2. Lemma}
If $\b$ is an $\epsilon$-number then $\club_\b$ is a proper ideal.
\endproclaim

\demo {Proof}
A good exercise for an undergraduate student. \qed
\enddemo

I will need to study the partial order $\CL_\b$ of the $\club_\b$-positive Borel subsets of $\C_\b$ ordered by inclusion. This forcing adds a closed subset of the ordinal $\beta$ of ordertype $\beta$ which has finite intersection with all ground model sets of ordertype $\o.$ 
I will show that this partial order is proper and its effect at cardinal invariants is limited to increasing the bounding number $\frak b$. 
Furthermore, as the ordinal $\b$ increases, these forcings form a tower that is a suitable candidate for adding a closed unbounded set of $\o_1$ that escapes guessing by the ground model sets.

First, some definitions due to Todorcevic. For every countable limit ordinal $\g$ fix an increasing function $c_\g\colon \o\to\g$ with a cofinal range. 
For countable ordinals $\xi\in\eta$ the {\bf walk} from $\eta$ to $\xi$ is a sequence $\eta=\eta_0\ni\eta_1\ni\dots\eta_k=\xi$ such that for every number $m\in k,$ if $\eta_m$ is a successor ordinal then $\eta_{m+1}$ is its predecessor, and if $\eta_m$ is a limit ordinal then $\eta_{m+1}$ is
the smallest element of $\rng(c_{\eta_m})$ greater or equal to $\xi.$ 
The {\bf norm} of such a
walk will be denoted by $n(\xi,\eta)$ and it is defined as the maximum of  the set
$\{n\colon $ for some $m\in k$ $\eta_m$ is a limit ordinal and $\eta_{m+1}=c_{\eta_m}(n)\}$. 
For every function $f\in\o^\o$ and every countable limit ordinal
$\g$ let $d(\g, f)$ be the set  of all $\xi\in\gamma$ such that $n(\xi,\eta)\in f(k)$ where $\eta=c_\g(k)$ is the
smallest element of the set $\rng(c_\g)$ above $\xi$.
The following is not difficult to check\colon 

\roster
\item $d(\g, f)$ is a subset of $\g$ of ordertype (less or) equal to $\o$
\item every cofinal subset of $\g$ of ordertype $\o$ is a subset of some $d(\g, f)$
\item if $f<g$ up to a finitely many values then $d(\g, f)\subset d(\g, g)$ up to finitely
many ordinals.
\endroster

Let me just indicate the proof of (1). 
It is enough to show that there are only finitely many elements of the set $d(\g, f)$
between $c_\g(k-1)$ and $c_\g(k).$ 
Look at the tree of all possible walks starting at $c_\g(k)$ with
norm less than $f(k).$ 
This is a finitely branching tree with no infinite branch, therefore it must be finite, and
all elements of the set $d(\g,f)\cap[c_\g(k-1),c_\g(k))$ are mentioned in this tree. 
Thus the set is finite and (1)
follows.

\proclaim {1.3. Lemma}
For every $\epsilon$-number $\b$, $\cov(\club_\b)=\add(\club_\b)=\frak b.$
\endproclaim

\demo {Proof}
To show that $\add(\club_\b)\geq \frak b$ it is enough to show that given $\kappa\in\frak b$ many subsets
$d_\alpha\colon \alpha\in\kappa$ of $\b$ of ordertype $\o,$ then the set $\bigcup_{\a\in\kappa}A_{d_\a}$
belongs to the ideal $\club_\b.$ 
For every ordinal $\a\in\k$ let $\g_\a=\sup(d_\a)$ and choose a function
$f_\a\in\o^\o$ such that $d_\a\subset d(\g_\a,f_\a).$ 
There is a function $g\in\o^\o$ modulo finite dominating all the functions $f_\a.$ 
Then for every $\a\in\k$ the set $d_\a$ is included in $d(\g_\a, g)$
up to a finite set, and so the set  $\bigcup_{\a\in\kappa}A_{d_\a}$ is a subset of the $\club_\b$-small
set $\bigcup_{\gamma\in\b+1}A_{d(\g, g)}$ as desired.

For $\cov(\club_\b)\leq\frak b$ just choose a modulo finite unbounded collection $f_\a\colon \a\in\frak b$
of increasing functions in $\o^\o$ and let $d_\a=d(\b,f_\a).$ 
I will argue that the sets $A_{d_\a}$
together cover the whole space $\C_\b.$ 
Let $e\subset\b$ be an arbitrary cofinal set. 
Let $x\subset\o$ be the set of those integers $k$ such that $e$ has some element between the ordinals $c_\b(k-1)$ inclusive and $c_\b(k)$, and let $g\colon x\to\o$ be defined as $g(k)=$ the norm of the walk from $c_\b(k)$ to the minimum of $e\setminus c_\b(k-1)$.
Since the functions $f_\a$ are increasing, there must be an ordinal $\a\in\frak b$ such that the set
$\{k\in x\colon g(k)\in f_\a(k)\}$ is infinite. By the definitions, $e\cap d(\b, f_\a)$ is infinite and
$e\in A_{d(\b,f_\a)}$ as desired. \qed
\enddemo

Digging deep into the present paper, namely into the proof of Lemma 3.4, the reader should observe that Theorems 0.2 and 0.3 are based on the following simple consequence of Lemma 1.3\colon  for every $\club_\b$-positive Borel set $B\subset\C_\b,$
$\cov(I\restriction B)=\cov(I)=\frak b.$  
It is exactly this homogeneity property, or the lack of thereof, that makes the Ostaszewski principle and others much harder to deal with than WCG.

\proclaim {1.4. Lemma}
$\frak b=\aleph_1$ implies the weak club guessing principle.
\endproclaim

\demo {Proof}
Suppose that $f_\a\colon \a\in\omega_1$ is a modulo finite unbounded collection of increasing functions in $\o^\o.$
The last paragraph of the previous proof immediately shows that the collection $d(\g, f_\a)\colon \a\in\o_1,\g\in\o_1$ limit, exemplifies the weak club guessing principle WCG. \qed
\enddemo

Typically, $\frak b=\aleph_1$ implies a statement much stronger than WCG. 
Actually all infinite subsets of $\omega_1$
will have infinite intersection with one set in the oracle presented in the previous proof. 
I do not know whether this behavior occurs in other prediction principles.

\proclaim {1.5. Lemma} 
For an $\epsilon$-number $\b$ the forcing $\CL_\b$ is proper.
\endproclaim

\demo {Proof}
Faced with several options, I will present a proof that provides literally no insight into the combinatorics of the forcing.

\proclaim {1.6. Claim}
Whenever $P$ is a forcing adding a closed subset $\dot e$ of $\b$ of ordertype $\b$
which has a finite intersection with every ground model set of ordertype $\o,$
and $M$ is a countable elementary submodel of a large enough structure containing
some condition $p\in P,$ the set $\{\dot e/G\colon G\subset M\cap P$ is an $M$-generic filter, $p\in G\}$ is $\club_\b$-positive.
\endproclaim

The lemma immediately follows by an absoluteness argument. 
Namely, if $p\in\CL_\b$
is a condition and $M$ is a countable elementary submodel of a large enough structure
containing the condition $p,$ the set $q=\{e\in p\colon  e$ is an $M$-generic set for $\CL_\b\}$ is a Borel set, it is $\club_\b$-positive set by the claim, and as such it
is a condition in the forcing $\CL_\b$. 
By Borel absoluteness, this condition forces
the generic set to be actually $M$-generic, so it is the required master condition for the
model $M.$

To prove the claim, fix a condition $p\in P$ and consider the following infinite game between players Adam and Eve.

\roster
\item Adam produces on a fixed schedule, finite piece by finite piece, subsets $d_n\colon n\in\omega$ of $\b$ of ordertype $\o,$ one by one ordinals $\g_n=\sup(d_n)$ and one by one open dense subsets $O_n\colon n\in\omega$ of the forcing $P.$ Literally this means that at round $n$ he indicates the ordinal $\g_n\in\b+1,$ the open dense set $O_n$ and the finite collection of finite sets $d_m^n\colon m\in n$ such that in the future $d_m\cap c_{\g_m}(n)=d_m^n,$ for all numbers $m\in n.$
\item Eve produces a descending chain $p_n\colon n\in\o$ of conditions below $p.$ Eve is allowed
to tread water, that is, to wait for any finite number of rounds before producing the next condition, but
she must make sure that $p_{n+1}\in O_n.$ 
\endroster

Eve wins if the filter $g\subset P$ generated by  the set of her answers meets all the dense sets necessary so that the expression $\dot e/g$ makes sense, and the set $\dot e/g$ has finite intersection with all the sets $d_n\colon n\in\o.$ 
The game is Borel and therefore determined by \cite {M}.

\proclaim {1.7. Claim}
Eve has a winning strategy.
\endproclaim

\demo {Proof}
If this is not the case then Adam must have a winning strategy $\s$. I will derive a contradiction.

First, a small observation. 
Given a condition $q\in P$ and a finite set $\g_m\colon m\in n$ of limit ordinals below
$\b+1,$ there must be a number $k$ such that for every collection $d_m\colon m\in n$ of sets of ordertype $\o$
with $\sup(d_m)=\g_m,$ there is a condition $\bar q\leq q$ forcing $\forall m\in n\ \dot e\cap\check d_m\subset
\check c_{\g_m}(k).$ 
If this failed for some condition $q$ and ordinals $\g_m\colon m\in n,$ for each number $k\in \o$ there would be a conterexample $d_m^k\colon m\in n.$ 
For each number $m\in n$ let then $b_m$ be the ``diagonal union" of the sets $d_m^k,$
so that $b_m=\bigcup_{k\in l} d_m^k$ as far as the ordinals between $c_{\g_m}(l-1)$ inclusive and $c_{\g_m}(l)$ are concerned, this for all numbers $l\in\o.$ 
By the properties of the name $\dot e$, there then must exist a condition
$\bar q\leq q$ and a number $k$ such that $\bar q\Vdash\forall m\in\check n\ \dot e\cap\check b_m\subset c_{\g_m}(\check k).$
This contradicts the choice of the sets $d_m^k\colon m\in n.$

Now let Eve face the strategy $\s.$ 
By induction on $n\in\o$ she will construct a counterplay with intermediate positions $0=\tau_0\subset\tau_1\subset\dots$ and she will also create a log of conditions $p=p_0\geq p_1\geq\dots$ and numbers
$k_0\in k_1\in\dots$ so that the following inductive conditions are satisfied\colon 

\roster
\item $\tau_{n+1}$ obtains from $\tau_n$ by Eve's waiting for some time and then playing $p_{n+1}.$
\item $p_{n+1}\in O_n$ and $p_{n+1}$ decides whether the $n$-th ordinal below $\b$ in some fixed enumeration belongs to the set $\dot e$ or not.
\item at the position $\tau_n$ the strategy $\sigma$ has already decided what the finite set $d_m\cap c_{\g_m}(k_n)$ is going to be, for each $m\in n.$
\item letting $d_m^n$ be the finite set the strategy $\sigma$ indicated as $d_m\cap c_{\g_m}(k_n),$ $p_{n+1}\Vdash\dot e\cap\check d_m^n\subset c_{\g_m}(k_{m+1}),$ this for all $m\in n.$
\endroster

If Eve succeeds in doing this, then in the end she won\colon  writing $g$ for the filter generated by her answers, (4) shows that for all numbers $n,$ $\dot e/g\cap d_n\subset\g_n(k_{n+1})$ and since the set $d_n$ is cofinal in $\g_n$ of ordertype $\o,$
the intersection $\dot e/g\cap d_n$ is finite, as desired. 
This will complete the proof of the Claim. 
But in order to maintain the inductive conditions (1--4) Eve also needs to maintain the following\colon 

\roster
\item "{(5)}" for every collection $b_m\colon m\in n$ of cofinal subsets of ordinals $\g_m$ of ordertype $\o$ there is a condition $q\leq p_n$ such that for each $m\in n$ $q\Vdash\dot e\cap\check b_m\subset c_{\g_m}(k_n).$
\endroster

To start the construction, Eve puts $\tau_0=0, p_0=p$ and $k_0=0.$ 
The items (1--5) are vacuously satisfied. Suppose $\tau_n, p_n, k_n$ have been constructed. 
Then Eve looks at the objects $d_m\colon m\in\o,\g_m\colon m\in\o, O_m\colon m\in\o$ the strategy $\s$ produces if she treads water from the position $\tau_n$ on.
 (``What would you do if you did not have me?") 
By (5), there is a condition $q\leq p_n$ such that for all $m\in n$ $q\Vdash\dot e\cap\check d_m\subset \check c_{\g_m}(k_n).$ Let $p_{n+1}\leq q$ be some condition in the dense set $O_n$ that decides whether the $n$th ordinal below $\beta$ in some fixed enumeration belongs to the set $\dot e$ or not. 
Let $k_{n+1}$ be a number above $k_n$ such that for all sets $b_m\colon m\in n+1$ cofinal in the respective ordinals $\g_m\colon m\in n+1$ of ordertype $\o$ there is a condition $\bar q\leq p_{n+1}$ forcing $\dot e\cap b_m\subset c_{\g_m}(k_{n+1})$ for all $m\in n+1.$ 
Such a number exists by the observation in the beginning of this proof. 
The position $\tau_{n+1}$ is now obtained from $\tau_n$ by Eve waiting until the strategy $\s$ commits on the finite sets
$d_m\cap c_{\g_m}(k_{n+1})\colon m\in n+1$ as well as on the ordinal $\g_n$ and the open dense set $O_n,$ and then plays $p_{n+1}.$ 
The inductive conditions (1--5) are satisfied and the Claim follows.
 (``Why, of course you would never get anywhere!") \qed
\enddemo

Claim 1.6 immediately follows from Claim 1.7. 
Suppose $p\in P$ is a condition and $\s$ is Eve's winning strategy in the game starting at $p.$ 
Whenever $M$ is a countable elementary submodel of large enough structure containing the strategy $\s,$ the set $\{\dot e/G\colon G\subset M\cap P$ is an $M$-generic filter, $p\in G\}$ is Borel by \cite {Z, Claim 1.1.3} and $\club_\b$-positive. For if it were small, it would be covered by countably many sets $A_{d_n}\colon n\in \o$ for some sets $d_n\subset\b$ of ordertype $\o.$ 
Then let Adam play against the strategy $\s,$ on the way creating these sets $d_n$ and enumerating all open dense subsets of the poset $P$ in the model $M.$ 
Since the strategy $\s$ is winning for Eve, her answers must have generated an $M$-generic filter  $G\subset P\cap M$ such that $\dot e/G\notin\bigcup_{n\in\o}A_{d_n}$ as desired. \qed
\enddemo

\proclaim {1.8. Lemma}
(ZFC+LC) Let $\b$ be an $\epsilon$-number. 
Every universally Baire $\club_\b$-positive subset of $\C_\b$ has a Borel positive subset.
\endproclaim

\demo {Proof}
I will present a proof which does not mention determinacy except for the reference to
Claim 1.6. 
I will show that given an inaccessible cardinal $\k,$ in the Solovay model obtained from $\k$
it is the case that every $\club_\b$-positive subset of $\C_\b$ has a Borel positive subset. 
This will be enough. 
It is well known that if $A$ is a universally Baire set and $T, S$ are class trees projecting
to complements in all set generic extensions such that $A=p[T],$ then the theory of the model $\langle L(\R)[p[T]],\in, p[T]\rangle$ is absolute throughout generic extensions satisfying DC, if suitable large cardinals exist \cite {W}. 
If such a set $A$ is positive, the theory of $L(\R)[p[T]]$ in the Solovay model must see that $p[T]$ is a positive set and that it has a Borel positive subset.
 Therefore the same must happen in the ground model, proving the lemma.

So suppose that $\k$ is an inaccessible cardinal, $G\subset Coll(\o,<\k)$ is a generic filter and $A\subset\C_\b$ is a positive set in the Solovay model $V(\R\cap V[G])\subset V[G].$ 
Thus there is a formula $\phi$, a real $r\in V[G]$ and an element $t\in V$
such that $V[G]\models A=\{e\in\C_\b\colon \phi(r,t,e)\}$. 
Since the set $A$ is positive, in the model $V[r]$ there must be a partial order
$P$ of size less than $\k$ and a $P$-name $\dot e$ for a closed
subset of $\b$ of ordertype $\b$ with finite intersection with all $V[r]$-sets of ordertype $\o,$ such that $P\Vdash Coll(\o,<\k)\Vdash\phi(\check r,\check t,\dot e).$

Now apply the proof of Claim 1.8 to see that the set $B=\{e\in\C_\b\colon e$ is $V[r]$-generic for the poset $P\}$
is Borel and $\club_\beta$-positive. 
Also $B\subset A$, completing the proof of the lemma. \qed

\enddemo

\head {2. The tower of the ideals $\club_\b$}\endhead

The forcing $\CL$ for adding a closed unbounded subset of $\o_1$ which
escapes guessing by the ground model sets is obtained as the tower of
the posets $\CL_\b$. 
The situation is somewhat analogous to the
nonstationary tower forcing.

\definition {2.1. Definition}
The forcing $\CL$ is the union of all forcings $\CL_\b$ for (countable) $\epsilon$-numbers
$\b$. 
For a condition $p\in\CL$ write $\b_p$ for the unique countable ordinal $\b$ such that $p\in \CL_\b.$ 
The ordering is defined by $q\leq p$ if $\b_p\leq\b_q$ and for every set
$e\in q,$ $e\cap\b_p\in p.$ 
Thus if $\b_p\in\b_q$ and $q\leq p$ then $\b\in\bigcap q.$
\enddefinition

In this section I will show that the $\CL$ extension is canonically given by a
closed unbounded set $E_{gen}\subset\o_1$ and that the forcing $\CL$ is proper.
In the next section I will argue that the countable support iteration of the
forcing $\CL$ isolates the cardinal invariant $\frak b$ and forces the failure
of the weak club guessing principle, proving Theorem 0.3.

\definition {2.2. Definition}
The $\CL$-name $\dot E_{gen}$ is defined as the name for the set of those ordinals $\g\in\check\o_1$ for which there is a condition $p$ in the
generic filter such that $\gamma\in\bigcap\dot p.$
\enddefinition

\proclaim {2.3. Lemma}
\roster
\item $\CL\Vdash$ the generic filter is the set of those conditions
$p\in\check\CL$ such that $\dot E_{gen}\cap\check\b_p\in\dot p$.
\item $\CL\Vdash\dot E_{gen}\subset\check\o_1$ is a closed unbounded set which has a finite intersection with all ground model sets of ordertype $\omega.$
\endroster
\endproclaim

\demo {Proof} 
For the proof of (1) first observe that if $p\in\CL$ is a countable union $p=\bigcup_np_n$ of Borel sets then $p$ forces that one of the sets $p_n$ appears in the generic filter. 
For let $q\leq p.$ Writing $q_n=\{e\in q\colon e\cap\b_p\in p_n\},$ it is clear that $q=\bigcup_n q_n$ and so one of the sets $q_n$ is $\club_{\b_q}$-positive. 
Then the condition $q_n\leq q\leq p$ forces $p_n$ into the generic filter--it is stronger than $p_n.$

This means that if $G\subset\CL$ is a generic filter and $\b$ is an $\epsilon$-number such that for some condition
$p\in\CL_\b,$ $p\in G,$ then the ultrafilter $G\cap\CL_\b$ respects countable disjunction from the ground model. 
It is well known that for any such ultrafilter on the algebra of Borel subsets of a Polish space there is a singleton $\{e\}$ which is the intersection of all the sets in the ultrafilter.
It is clear that $e=E_{gen}\cap\b$--otherwise there would be an ordinal $\g\in\b$ such that $\{f\in\C_\b:\g\in f\}\in G\not\leftrightarrow\g\in E_{gen},$ contradicting the definition of the name $\dot E_{gen}.$ 
(1) follows.

(2) is left to the reader. \qed  
\enddemo

\proclaim {2.4. Lemma}
The forcing $\CL$ is proper.
\endproclaim

\demo {Proof}
The proof follows the lines of the argument for Lemma 1.5. 
Fix a condition $p\in\CL$ and consider the following infinite game between players Adam and Eve.

\roster
\item Adam produces on a fixed schedule, finite piece by finite piece, subsets $d_n\colon n\in\omega$ of $\o_1$ of ordertype $\o,$ one by one ordinals $\g_n=\sup d_n$ and one by one open dense subsets $O_n\colon n\in\omega$ of the forcing $\CL.$ Literally this means that at round $n$ he indicates the ordinal $\g_n,$ the open dense set $O_n$ and the finite collection of finite sets $d_m^n\colon m\in n$ such that in the future $d_m\cap\g_m(n)=d_m^n.$
\item Eve produces a descending chain $p_n\colon n\in\o$ of conditions below $p.$ 
Eve is allowed
to tread water, that is, to wait for any finite number of rounds before producing the next condition, but
she must make sure that $p_{n+1}\in O_n$ and $\g_n\in\b_{p_{n+1}}$.
\item Additionally, when Eve played $p_n,$ in the very next round Adam can play a finite set $x_n$ of countable ordinals larger than $\b_{p_n}.$ 
When Eve plays the condition $p_{n+1}$ she is required to do it so that $x_n\subset\b_{p_{n+1}}$
and that every set $e\in p_{n+1}$ is disjoint from the set $x_n.$
\endroster

Eve wins if the filter $g\subset\CL$ generated by  the set of her answers meets all the dense sets necessary so that the expression $\dot E_{gen}\cap\sup_n\g_n/g$ makes sense, and the set $\dot E_{gen}\cap\sup_n\g_n/g$ has finite intersection with all the sets $d_n\colon n\in\o.$ 
The game is Borel and therefore determined by \cite {M}. 
The proof of the following claim repeats the argument for Claim 1.7 almost verbatim.

\proclaim {2.5. Claim}
Eve has a winning strategy.
\endproclaim

Now let $p\in\CL$ be an arbitrary condition and let $M$ be a countable elementary submodel of a large enough structure containing the condition $p$ as well as some Eve's winning strategy $\sigma$ in the above game. Write $\b=M\cap\o_1.$ 
I will prove that the Borel set $q=\{e\in\C_\b\colon  e\cap\b_p\in p$ and $e$ is an $M$-generic club for the poset $\CL\}$ is $\club_\b$-positive. 
This will complete the proof of the lemma, since then $q\leq p$ will be a condition which by Borel absoluteness and Lemma 2.3(1) forces $\dot E_{gen}\cap\check\beta$ to be $M$-generic, therefore it is the required master condition for the model $M.$

So let $d_n\colon n\in\o$ be a countable collection of subset of $\b$ of ordertype $\o$; I must produce a set $e\in q$ which has finite intersection with all of them. 
Just as in the proof of Lemma 1.5, this set will obtain as a result of a suitable play against the strategy $\s$, however here a little more sophistication is needed to obtain the required play.

First of all, by rearranging the collection of the sets $d_n\colon n\in\o$ and taking diagonal unions if necessary, I can assume that it comes in the form $d_\g\colon \g\in\b+1$ limit so that $\sup(d_\g)=\g.$ 
Inside the model $M$ it is possible to find a suitable large structure $H_\theta$ and countable elementary submodels $N_0\in N_1\in\dots$ of it so that $p,\s\in N_0$ and $M\cap H_\theta=\bigcup_nN_n$. 
Of course, the $N$-sequence itself will not be in the model $M.$ 
And it is possible to find enumerations $\bar\g_n\colon n\in\o$ of countable limit ordinals in $M$ and $\bar O_n\colon n\in\o$ of open dense subsets of the poset $\CL$ in $M$ so that $\g_n$ and $\bar O_n$ are both in the model $N_n.$

By induction on $n\in\o$ build partial plays $\tau_0\subset\tau_1\subset\dots$ against the strategy $\s$ so that

\roster
\item $\tau_n\in N_n$ and the play $\tau_n$ ends with Eve finally making the move $p_n$
\item in the first round $m_n$ of the play $\tau_{n+1}$ after $\tau_n$ ended, Adam puts $x_n=N_{n+1}\cap d_\b\setminus N_n,$ $\g_{m_n}=\bar\g_n$ and $O_{m_n}=\bar O_n$
\item Adam plays so that in the end the infinite subset of $\bar\g_n$ he produced is exactly $d_{\bar\g_n}$.
\endroster

This is easily possible. 
Look at the filter $g$ generated by the conditions $p_n\colon n\in\o$ produced by the strategy $\s$ in the course of the play $\bigcup_n\tau_n.$ 
Clearly the filter $g\subset\CL$ is $M$-generic, contains the condition $p$ and $\dot E_{gen}\cap\b/g\cap d_\g$ is finite for all limit ordinals $\g\in\b.$ The last thing to verify is that the intersection $\dot E_{gen}\cap\b/g\cap  d_\b$ is finite. But Adam's moves $x_n$ were chosen exactly so as to guarantee that $\dot E_{gen}\cap \b/g\cap d_\b\subset N_0,$ and the model $N_0$ contains only a finite piece of the set $d_\b.$ \qed
\enddemo

\head {3. The iteration}\endhead

In this section I will finally present the proof of Theorem 0.3. 
Argue in the theory ZFC+LC. Look at the countable support iteration $P$ of the poset $\CL$
of length $\frak c^+.$ 
This iteration adds generic clubs $\dot E_\alpha\colon \alpha\in\frak c^+.$
By a standard properness and chain condition arguments, $P$ collapses $\frak c$ to $\aleph_1$ and leaves other cardinals standing, it
forces WCG to fail and so it makes $\frak b=\aleph_2$. 
Suppose that $\frak x$ is a tame cardinal invariant, defined as $\min\{|A|\colon A\subset\R,\phi(A)\land\forall x\in\R\exists y\in A\ \theta(x,y)\}$,
where the quantifiers in the formula $\phi$ range over natural numbers or elements of the set $A$ and $\theta$ is a projective formula making no mention of the set $A$, and suppose that $\frak x<\frak b$ holds in some forcing extension. 
I will prove that $P$ forces $\frak x=\aleph_1$ to hold. 
This will conclude the proof of Theorem 0.3.

Without loss of generality I can assume that the continuum hypothesis holds, because the forcing with
the first $\omega_1$ many copies of $\CL$ restores the continuum hypothesis and the whole
argument below can then be repeated in the resulting model.

First, let $\beta$ be an $\epsilon$-number and $\alpha$ be a countable ordinal. 
Arguing as in \cite {Z} using the results of Section 1 of the present paper, define a presentation of the countable support iteration of the posets $\CL_\b$ of length $\alpha$ with the following two general definitions and a lemma\colon

\definition {3.1. Definition}
Let $\C_\b^\a$ denote the collection of all $\a$-sequences of elements of $\C_\b$. The set $\CL_\b^\a$ consists of those nonempty Borel sets $p\subset\C_\beta^\alpha$ satisfying the following conditions\colon 

\roster
\item for every ordinal $\g\in\b$ the set $p\restriction\g=\{\vec e\restriction\g\colon \vec e\in p\}$ is Borel
\item for every ordinal $\g\in\b$ and every sequence $\vec e\in p\restriction\g$ the set of all sets $f\in\C_\b$
with $\vec e^\smallfrown f\in p\restriction\g+1$ is $\club_\b$-positive
\item whenever $\vec e_0\in p\restriction\g_0, \vec e_1\in p\restriction\g_1, \vec e_2\dots$ are sequences such that $\vec e_0\subset\vec e_1\subset\dots$ then $\bigcup_n\vec e_n\in p\restriction \bigcup_n\g_n.$  
\endroster

The ordering is that of inclusion.

\enddefinition

\definition {3.2. Definition}
The ideal $\club_\beta^\alpha$ on $\C_\b^\a$  is the collection of those sets $B\subset\C_\beta^\alpha$
for which Adam has a winning strategy in the game $G(B).$ 
The game $G(B)$
lasts $\alpha$ rounds, and at each round $\gamma\in\alpha$ Adam plays a Borel set $X_\gamma$ in the ideal $\club_\beta$ and subsequently Eve plays a set $e_\gamma\in\C_\beta$ that does not belong to the set $X_\gamma.$ 
Eve wins if the sequence $\langle e_\gamma\colon \g\in\a\rangle$ belongs to the set $B.$
\enddefinition

\proclaim {3.3. Lemma}
\roster
\item The poset $\CL_\b^\a$ is forcing equivalent to the countable support iteration of the posets
$\CL_\b$ of length $\a.$
\item (ZFC+LC) For every projective set $B\subset\C_\b^\a,$ either the set is $\club_\b^\a$-small
or there is a condition $q\in\CL_\b^\a$ with $q\subset B.$
\endroster
\endproclaim

The proof of the lemma follows the lines of \cite {Z, Section 1} verbatim, using Lemmas 1.5 and 1.8 of the present paper. 
Now back to our particular setup. 
The key step in the proof of Theorem 0.3 is

\proclaim {3.4. Lemma}
There is a set $A\subset\R$ such that $\phi(A)$ holds and for every $\epsilon$-number $\b$, every countable ordinal $\a$, every condition $p\in\CL_\b^\a$ and every Borel function $f\colon p\to\R$ there is a condition $q\leq p$ and a real $y\in A$ such that for every $x\in q$ $\theta(f(x),y)$ holds.
\endproclaim

\demo {Proof}
This is similar to \cite {Z, Section 2}. 
Note that the assertion to be proved is $\Sigma^2_1$,
so it is enough to verify it in some generic extension by the $\Sigma^2_1$-absoluteness
theorem of Woodin. 
Well, move into the postulated generic extension $V[G]$ satisfying
$\frak x<\frak b$ and there choose a set $A\subset\R$ of size $\frak x<\frak b$ such that
$\phi(A)$ and $\forall x\in\R\exists y\in A\ \theta(x,y)$ both hold. 
Working in the model
$V[G]$ I will show that the set $A$ has the required properties.

Let $\b,\a,p\in\CL_\b^\a,f\colon p\to\R$ be arbitrary as in the statement of the lemma. 
For every
real $y\in A$ let $B_y=\{\vec e\in p\colon \theta(f(\vec e),y)\}.$ 
If one of these sets is
$\club_\b^\a$-positive then by the dichotomy of Lemma 3.3(2) applied in the model $V[G]$ we are done.
So it is enough to derive a contradiction from the assumption that these sets are all
$\club_\b^\a$-small. 
In such a case, choose winning strategies $\sigma_y\colon y\in A$
for Adam in the respective games $G(B_y).$ 
By induction on $\g\in\a$ choose
sets $e_\g$ such that at every stage $\d\in\a+1$ of the induction, the sequence
$\langle e_\g\colon \g\in\d\rangle$ belongs to the set $p\restriction\d$ and is a
legal partial counterplay against all the strategies $\s_y.$ 
The inductive assumption
clearly persists at limit stages due to the countable support condition. 
To get
the next real $r_\d,$ observe that all the sets $\s_y(\langle e_\g\colon \g\in\d\rangle)\colon y\in A$
the various strategies $\s_y$ advise Adam to play are $\club_\b^\a$-small,
and there are only $\frak x<\frak b$ many of them. 
By Lemma 1.3, their union
is still small and does not cover the positive set $\{e\colon  \langle e_\g\colon \g\in\d\rangle^\smallfrown e\in p\restriction\d+1\}.$ Just choose $e_\d$ in this set and outside of the union.

In the end, look at the sequence $\vec e=\langle e_\g\colon \g\in\a\rangle$ and the real $x=f(\vec e).$ 
Since the sequence $\vec e$ is a legal counterplay against
all of the winning strategies $\sigma_y\colon y\in A,$ it must be the case that for every $y\in A,$
$\theta(x,y)$ fails. But this contradicts the assumed properties of the set $A.$ \qed
\enddemo

Let $A\subset\R$ be as in the previous Lemma. 
I will show $P\Vdash\phi(\check A)\land \forall x\in\R\exists y\in\check A\ \theta(x,y).$ This will complete the proof of Theorem 0.3; $\frak x=\aleph_1$ will hold in the $P$-extension too as witnessed by the set $A$ of size $\frak c^V=\aleph_1.$

Now of course $P\Vdash\phi(\check A)$ because of the low syntactical complexity
of the formula $\phi.$ 
However, in order to prove $P\Vdash\forall x\in\R\exists
y\in\check A\ \theta(x,y)$, I need to understand all the possible $P$-names
for reals. 
The key idea is the approximation of such names by $\CL_\beta^\alpha$
names for suitable ordinals $\beta$ and $\alpha$. 
This approximation
is facilitated by another general lemma\colon

\proclaim {3.5. Lemma}
\roster
\item Let $\b$ be an $\epsilon$-number, let $\a$ be a countable ordinal and let $\pi\colon \a\to\o_2$ be an order-preserving map. Then $\pi$ naturally extends to a map $\bar\pi\colon \C_\b^\a\to \C_\b^{\rng(\pi)}$ and to an order-preserving map $\bar\pi\colon \CL_\b^\a\to P$ such that for every condition
$p,$ $\bar\pi(p)\Vdash \langle \dot E_{\pi(\gamma)}\cap\check\beta\colon \gamma\in\check\alpha\rangle\in\bar\pi''\dot p.$
\item Let $r\in P$ be an arbitrary condition and let $M$ be a countable
elementary submodel of some large structure containing the condition $r$.
Writing $\beta=M\cap\omega_1$, $\alpha=o.t. M\cap\o_2$ and $\pi\colon \a\to M\cap\o_2$
for the unique order-preserving bijection, there is a condition $p\in\CL_\beta^\alpha$ such that for every sequence $\vec e\in p$ the sequence $\bar\pi(\vec e)$ is $M$-generic for the poset $P,$ compatible with the condition $r.$ 
Necessarily $\bar\pi(p)$ is a master condition for the model $M$ stronger than $r.$
\endroster
\endproclaim

\demo {Proof}
(1) should not need an argument. 
For a sequence $\vec e\in\C_\b^\a,$ $\bar\pi(\vec e)$ is just a reindexing of it. 
And for a condition $p\in\CL_\b^a,$ $\bar\pi(p)$ is the condition $r$ in the poset $P$ with domain $\pi''\a$ so that for each ordinal $\g,$ $r\restriction\pi(\g)\Vdash_{P_{\pi(\g)}}r(\pi(\g))=\{f\in\C_{\check \b}\colon \langle \dot E_{\check\pi(\d)}\colon \d\in\check\g\rangle^\smallfrown f\in\dot p\restriction\check\gamma+1\}$.

(2) is a completely standard countable support iteration argument using the proof of Lemma 2.4. 
Set $\pi(\a)=\o_2.$ 
By induction on $\g\in\a+1$ prove that\colon  (IH) for every ordinal $\d\in\g$, every condition $p_\d\in\CL_\b^\d$ such that for every sequence $\vec e\in p_\d$ the sequence $\bar\pi(\vec e)$ is $M$-generic for the poset $P_{\pi(\d)}$, and every condition $r\in M\cap P_{\pi(\g)}$ with $\bar\pi(p_\d)\leq r\restriction\pi_\d$ there is a condition $p_\g \in\CL_\b^\g$ such that
$p_\d=p_\g\restriction\d,$ $\bar\pi(p_\g)\leq r$ and for every sequence $\vec e\in p_\g$ is $M$-generic for the poset $P_{\pi_\g}.$ 
(2) is then the application of this general fact to $\d=0$ and $\g=\a.$ 
The last sentence of (2) follows by a borel absoluteness argument and (1).

So suppose that the induction hypothesis IH has been verified up to an ordinal $\g.$ 
There are two cases. 
First let $\g=\bar\g+1$ be a successor ordinal, and $\d\in\g, p_\d\in\CL_\b^\d$ and $r\in P_{\pi(\g)}$ be as in IH. Then use IH at $\bar\g$ to get the suitable condition $p_{\bar\g}\in\CL_b^{\bar\g}$ and let $p_\g=\{\vec e\in C_\b^\g\colon \vec e\restriction\bar\g\in p_{\bar\g}$ and $\vec e(\bar\g)$ is an $M[\bar\pi(\vec e\restriction\bar\g)]$-generic club for the poset compatible with the condition $r(\bar\g)/\bar\pi(\vec e\restriction\bar\g)\}$. 
It is easy modulo the proof of Lemma 2.4 to check that this is the required condition. 
Second, let $\g=\sup_n\d_n$ be a limit ordinal, a supremum of an increasing sequence of smaller ordinals, and let $\d_0=\d, p_\d\in\CL_\b^\d$ and $r\in P_{\pi(\g)}$ be as in IH. 
Disregarding the condition $p_\d$ for a second, enumerate the open dense subsets of the poset $P_{\pi(\g)}$ in the model $M$ by $O_n\colon n\in\o$ and by induction on $n$ obtain a sequence of conditions $r=r_0\geq r_1\geq r_2\geq\dots$ in $P_{\pi(\g)}\cap M$ such that $r_{n+1}\restriction\d_n=r_n\restriction\d_n$ and $r_n\restriction\d_n\Vdash_{P_{\d_n}}$ for some condition $s\in\check O_n,$ $s\restriction\d_n$ belongs to the generic filter on $P_{\d_n}$ and $s$ is equal to $r_{n+1}$ above $\d_n$. 
Now using IH at the ordinals $\d_n$ in their turn, find conditions $p_{\d_n}\in\CL_\b^{\d_n}$ as in IH so that $\bar\pi(p_{\d_n})\leq r_n\restriction\d_n$ and $p_{\d_n}=p_{\d_{n+1}}\restriction\d_n.$ 
In the end, $p_\g=\{\vec e\colon \forall n\in\o\ \vec e\restriction\d_n\in p_{\d_n}\}$ is the desired condition in $\CL_\b^\g$. \qed
\enddemo

Now suppose that $r_0\in P$ is a condition and $r_0\Vdash\dot x$ is a real number.
I will produce a stronger condition $r_1$ and a real $y\in A$ such that
$r_1\Vdash\theta(\dot x,\check y).$ 
This will complete the proof of Theorem 0.3.
Well, choose a countable elementary submodel $M$ of a large enough structure
containing the condition $r_0$ and the name $\dot x.$ 
Write $\beta=M\cap\o_1,$ $\a=o.t.M\cap\o_2$ and $\pi\colon \a\to M\cap\o_2$ for the unique order-preserving
bijection and find a condition $p\in\CL_\b^\a$ as in Lemma 3.5(2). 
Let $f\colon p\to\R$ be the Borel function defined by $f(\vec e)=\dot x/\bar\pi(\vec e).$ Note that the latter expression makes sense since the sequence $\bar\pi(\vec e)$ is $M$-generic for the poset $P.$ 
By the properties of the set $A,$ there must be a condition $q\leq p$ in the forcing $\CL_\b^\a$ and a real $y\in A$ such that for every sequence $\vec e\in q$ $\theta(f(\vec e), y)$ holds. 
Consider the
condition $r_1=\bar\pi(q)\leq\bar\pi(p)\leq r_0$ in the forcing $P.$

It follows from the choice of $q,$ Lemma 3.5(1) and a Borel absoluteness argument that $r_1\Vdash\dot x=\dot f(\langle \dot E_{\pi(\gamma)}\cap\check\beta\colon \gamma\in\check\alpha\rangle)$ and $\theta(\dot x,\check y).$ 
This completes the proof of Theorem 0.3.

Theorem 0.2 immediately follows. 
If $\frak x$ is a tame cardinal invariant
such that ZFC+LC does not prove $\frak b\leq\frak x,$ then in the consistent
theory ZFC+LC+$\frak x<\frak b$ one can apply Theorem 0.3 to see that some partial order forces ZFC+LC+$\frak x=\aleph_1$+WCG fails. 
Thus ZFC+LC does not prove $\frak x=\aleph_1\to$WCG.

\enddocument